\documentclass[leqno,12pt]{amsart} 
 \usepackage{amssymb,amsmath} 
\theoremstyle{plain} 
 \newtheorem{theorem}{\indent\sc Theorem}[section]
 \newtheorem{lemma}[theorem]{\indent\sc Lemma}

\theoremstyle{definition} 
 \newtheorem{definition}[theorem]{\indent\sc Definition}

\usepackage{pgf,tikz,mathrsfs}
\usepackage{graphicx,subfigure}
\usetikzlibrary{arrows}
\usepackage{soul}

\begin{document}

\title{New Characterizations of the Helicoid in a Cylinder}

\author{Eunjoo Lee}


\email{eunjoolee@ssu.ac.kr}

\subjclass[2010]{ 
Primary 53A10; Secondary 49Q05.
 }
 %
 \keywords{ 
Helicoid, Rigidity}

 \thanks{ 
$^{*}$The author would like to express her deep gratitude to Professor Jaigyoung Choe for valuable comments and guidance. 
 }

\maketitle

\begin{abstract}
This paper characterizes a compact piece of the helicoid $H_C$ in a solid cylinder $C \subset \mathbb{R}^3$ from the following two perspectives. First, under reasonable conditions, $H_C$ has the smallest area among all immersed surfaces $\Sigma$ with $\partial \Sigma \subset d_1 \cup d_2 \cup S$, where $d_1$ and $d_2$ are the diameters of the top and bottom disks of $C$ and $S$ is the side surface of $C$. Second, other than $H_C$, there exists no minimal surface whose boundary consists of $d_1$, $d_2$, and a pair of  \textcolor{black}{rotationally symmetric} curves $\gamma_1$, $\gamma_2$ on $S$ along which it meets $S$ orthogonally. We draw the same conclusion when the boundary curves on $S$ are a pair of helices of a certain height. 
\end{abstract}

\section*{Introduction}

The helicoid in $\mathbb{R}^3$ is a classic example of a minimal surface. It is a simply connected, complete, embedded, ruled minimal surface foliated by helices and having infinite total curvature. It is symmetric with respect to its central axis and any horizontal line it contains. Catalan \cite{Catalan:1842} verified that, other than the plane, the helicoid is the unique ruled minimal surface. Collin and Krust \cite{Collin_Krust:1991} showed that a nonplanar minimal surface bounded by two lines, whose interior is a graph over a band, must be part of the helicoid. Colding and Minicozzi \cite{Colding_Minicozzi:2004} proved that every embedded minimal disk in $\mathbb{R}^3$ is either a graph of a function or part of an appropriately scaled helicoid. Based on their work, Meeks and Rosenberg \cite{Meeks_Rosenberg:2005} showed that a complete, embedded, simply connected nonplanar minimal surface must be the helicoid. Recently, Choe and Hoppe \cite{Choe_Hoppe:2013} constructed a higher-dimensional helicoid, whose further generalization is demonstrated in \cite{Lee_Lee:2017}.\

Bernstein and Breiner \cite{Bernstein_Breiner:2014} proved that part of the catenoid has the smallest area among the embedded minimal annuli in a slab spanning two parallel planes in $\mathbb{R}^3$. Inspired by their work, Choe conjectured that part of the helicoid minimizes the area among the surfaces spanning two skew diameters of the top and bottom disks of a solid cylinder $C$, with other boundaries lying on the side of $C$. His question was natural, because the catenoid is locally isometric to the helicoid and, therefore, if a compact piece of the catenoid has an area-minimizing property in a certain setting, then a compact piece of the helicoid is highly likely to inherit a similar property correspondingly. \

Because of the existence of trivial counterexamples, however, we had to impose additional (yet reasonable) conditions to obtain meaningful results. We were able to prove two classical results on the least area property of a piece of helicoid and three theorems on the uniqueness of the helicoid in a certain setting. We should point out that our results concern some properties of a compact portion of the helicoid, whereas previous characterization results were mostly about the complete helicoid. Additionally, our proofs use classical methods, and are very geometrical. 

The main results of this paper are twofold: concerning the least area property and the uniqueness of a piece of the helicoid centered in a cylinder in $\mathbb{R}^3$. First, we investigate the conditions under which a piece of the helicoid $H_C$ in a right circular compact cylinder $C$ has the smallest area among the surfaces spanning two skew diameters $d_1$ and $d_2$ of the top and bottom disks and bounded by the curves on the side surface $S$ of $C$. Second, we prove that $H_C$ is the unique minimal surface spanning $d_1$ and $d_2$ and meeting $S$ along the  \textcolor{black}{rotationally symmetric} curves. Furthermore, we determine other conditions that guarantee the uniqueness of $H_C$. \

Note that, in this paper, we consider every surface to be two-dimensional, orientable, and immersed in the three-dimensional Euclidean space $\mathbb{R}^3$.

\section{Helicoid as the least area surface}

\begin{figure}[h]
\centering
\subfigure[$H_C$]{\includegraphics[width=0.3\textwidth]{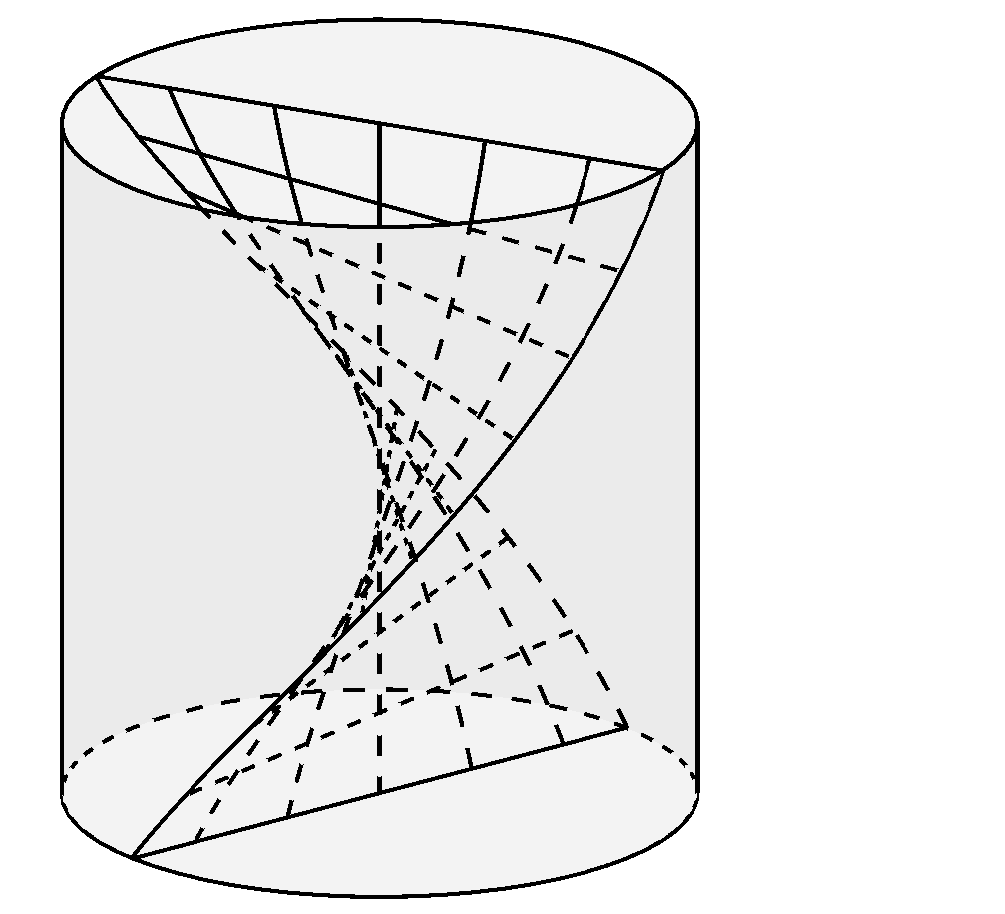}} 
\qquad \qquad
\subfigure[Two half disks with a neck]{\includegraphics[width=0.3\textwidth]{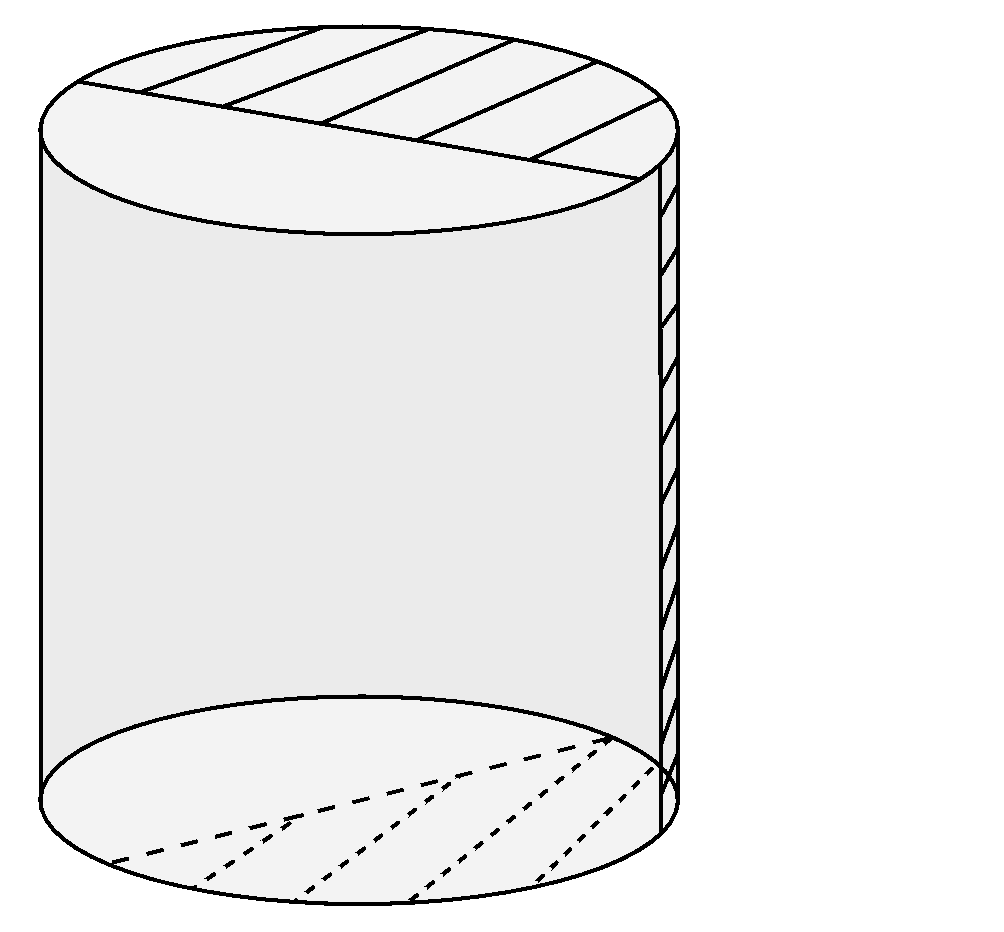}} 
\caption{Candidates for the Area-minimizer}
\end{figure}

Let $C \subset \mathbb{R}^3 $ be a right circular solid cylinder bounded by two disks $D_1$ and $D_2$ and a cylindrical surface $S$. Let $H$ be the helicoid ruled by lines perpendicular to the axis of $C$. We define the compact helicoid as $H_C=H \cap C$, whose boundary consists of two diameters $d_1 \in D_1$ and $d_2 \in D_2$ and two symmetrical helices $h_1$, $h_2$ $ \in S$. In fact, there are infinitely many such helicoids, but $H_C$ is defined as  \textcolor{black}{the one whose boundary helices spiraling with an angle less than or equal to $\pi/2$}. (See Figure 1(a)). Then, as a minimal surface, $H_C$ is a strong candidate for the least area surface among the immersed surfaces passing through $d_1$ and $d_2$ who have their boundaries on $S$,  \textcolor{black}{if $d_1$ and $d_2$ are not parallel to each other. It is clear that the plane passing through $d_1$ and $d_2$ is the area-minimizer when $d_1$ and $d_2$ are parallel. Assume, therefore, that $d_1$ and $d_2$ are not parallel.}

Depending on the height $h$ and radius $r$ of $C$, however, the surface composed of two half disks in $D_1$ and $D_2$ connected by a band of length $\epsilon>0$ in $S$ (Figure 1(b)) may have a smaller area than $H_C$ if $\epsilon$ is sufficiently small. In fact, when $h=r=1$ and the angle between $d_1$ and the vertical projection of $d_2$ onto $D_2$ is $\frac{\pi}{2}$, the area of the surface in Figure 1(b) is $\pi+\epsilon$, whereas the area of $H_C$ is $\pi(\sqrt{2}+\log(1+\sqrt{2}))>\pi+\epsilon$ for arbitrary small $\epsilon>0$. Therefore, it is natural and meaningful to seek the conditions under which $H_C$ has the smallest area. The next theorem provides a geometric condition on the candidate surfaces to rule out the trivial counterexample. By a disk-type surface we mean a surface topologically equivalent to a disk.

\begin{theorem} \label{thm1}
Let $\Sigma \subset C$ be an immersed disk-type surface with $\partial \Sigma \subset d_1 \cup d_2 \cup S$.  \textcolor{black}{Let $C_t=\{(x,y,z) \mid x^2+y^2 \leq t^2\}$.} If $\Sigma$ contains the axis of $C$  \textcolor{black}{and $\Sigma$ is transversal to $\partial C_t$ for every $0< t \leq r$}, then $\mathcal{H}^2 (\Sigma) \geq \mathcal{H}^2 (H_C)$. Equality holds if and only if $\Sigma$ is congruent to $H_C$.
\end{theorem}

To prove the theorem, we need the following lemma, which relates an integral of a function over a surface to the integral over the level sets. 

\begin{lemma} [Co-area formula] If $\Sigma$ is a  \textcolor{black}{Riemannian} manifold and $g: \Sigma \rightarrow \mathbb{R}$ is a proper (i.e. $g^{-1}((-\infty, t])$ is compact for all $t \in \mathbb{R}$) Lipschitz function on $\Sigma$, then, for any locally integrable function $f$ on $\Sigma$ and $t \in \mathbb{R}$, 

\[
\int_{\{h \leq t\}} f \, |\nabla_{\Sigma} \, g| \,dV = \int_{-\infty}^{t} \left(\int_{g=\tau} f \,dA_{\tau}\right) d\tau, 
\]
 \textcolor{black}{where $\nabla_{\Sigma}$ is the surface gradient on $\Sigma$.}
\end{lemma}

\begin{proof} \textcolor{black}{Take $N=\mathbb{R}$ in the co-area formula presented in \cite{Chavel:2001}.}
\end{proof}

\begin{proof}[Proof of Theorem \textcolor{black}{1.1}]
Without loss of generality, we may take $C=\{(x,y,z) \mid x^2+y^2 \leq r^2, \ 0 \leq z \leq h\}$ \textcolor{black}{and $H_C$ is centered at $z$-axis.}  \textcolor{black}{And let $l$ be the axis of $C$.  Then, for $d(x,y,z)=\sqrt{x^2+y^2}$, the transversality condition implies that $|\nabla_{\Sigma}d| \neq 0$ on $\Sigma \setminus l$. Therefore, $\displaystyle f:=\frac{1}{|\nabla_{\Sigma}d|}$ is locally integrable on $\Sigma \setminus l$. \textcolor{black}{Let $\hat{\Sigma}=\Sigma\setminus C_{\epsilon}$ for an arbitrarily small $\epsilon >0$.} For $A_{\textcolor{black}{\hat{\Sigma}}}(t)=\mathcal{H}^2(\textcolor{black}{\hat{\Sigma}} \cap C_t)$ \textcolor{black}{apply co-area formula to have}}

\[
A_{\textcolor{black}{\hat{\Sigma}}}(t)=\int_{\textcolor{black}{\hat{\Sigma}} \cap C_t} \, dS= \int_{\textcolor{black}{\epsilon}}^t \int_{\textcolor{black}{\hat{\Sigma}} \cap \partial C_\tau} \frac{1}{|\nabla_{\textcolor{black}{\hat{\Sigma}}}d|} \, ds \, d\tau.
\]
 
\begin{equation*}
\therefore \frac{d}{dt}A_{\textcolor{black}{\hat{\Sigma}}}(t)=\int_{\textcolor{black}{\hat{\Sigma}} \cap \partial C_t} \frac{1}{|\nabla_{\textcolor{black}{\hat{\Sigma}}}d|}\, ds \geq  \mathcal{H}^1(\textcolor{black}{\hat{\Sigma}} \cap \partial C_t),
\end{equation*}

since $|\nabla_{\hat{\Sigma}} d| \leq 1$. 

\begin{equation*}
\therefore \frac{d}{dt}A_{\textcolor{black}{\hat{\Sigma}}}(t) \geq \mathcal{H}^1(\textcolor{black}{\hat{\Sigma}}\cap \partial C_t).
\end{equation*}

$\it Claim.$ If $\Sigma$ contains the axis $l$ of $C$, then $\mathcal{H}^1(\Sigma \cap \partial C_t) \geq  \mathcal{H}^1(H_C \cap \partial C_t)$ for each $t \in(0,r]$. \\

To verify \textcolor{black}{the} $Claim$, let \textcolor{black}{$P_{\tau}=\{(x, y, z) \mid z=\tau\} \, (0 \leq \tau \leq h)$} be the plane parallel to $D_1$ and $D_2$. Given that $\Sigma$ is of disk-type and contains the axis of $C$, $\Sigma \cap P_{\tau}$ should contain at least one curve joining two points on the circle $\{(x,y,z)\mid x^2+y^2=r^2, \ z=\tau \}$ and passing through the center. In other words, if we slice $\Sigma$ by the planes parallel to $D_i$ (i=1,2), then at each height at least one of the intersection curves should join two points on the boundary circle and pass through the center. Because this is true for any $\tau \in [0, h]$, if we slice $C$ into the cylinders $C_t$ with radius $t \in (0,r]$, then $\Sigma \cap \partial C_t$ must have curves $\gamma_{p,t}$ joining $p^1_t$, $p^2_t$ and $\gamma_{q,t}$ joining $q^1_t$, $q^2_t$, where $d_i \cap C_t=\{p^i_t, q^i_t\}$, $i=1,2$. It is therefore clear that $ \mathcal{H}^1 (\Sigma \cap \partial C_t) \geq \mathcal{H}^1(\gamma_{p,t})+\mathcal{H}^1(\gamma_{q,t}) \geq \mathcal{H}^1(\textcolor{black}{h_{p,t}})+\mathcal{H}^1(\textcolor{black}{h_{q,t}})$ \textcolor{black}{where $h_{p,t}$ and $h_{q,t}$ are helices joining $p^1_t$, $p^2_t$ and $q^1_t$, $q^2_t$, respectively, because} the helix gives the shortest path between two points on the cylindrical surface (other than the vertical line). This proves the claim. \\

Thus, we have 
\[
\frac{d}{dt}A_{\textcolor{black}{\hat{\Sigma}}}(t) \geq \mathcal{H}^1 (\textcolor{black}{\hat{H}}_C \cap \partial C_t)=\int_{\textcolor{black}{\hat{H}_C} \cap \partial C_t}ds = \frac{d}{dt}A_{\textcolor{black}{\hat{H}}_C}(t).
\]
where $\hat{H}_C=H_C \setminus C_{\epsilon}$.\\
Integrating both sides over [$\epsilon, r$] yields $\textcolor{black}{A_{\hat{\Sigma}}(r) -A_{\hat{\Sigma}}(\epsilon)}\geq A_{\textcolor{black}{\hat{H}}_C}(r)\textcolor{black}{-A_{\hat{H}_C}(\epsilon)}$. Therefore, letting $\epsilon \to 0$ on both sides, we have \textcolor{black}{$A_{\Sigma}(r) \geq A_{H_C}(r)$, in other words,} $\mathcal{H}^2(\Sigma) \geq \mathcal{H}^2 (H_C)$. Note that the equality holds if and only if $|\nabla_{\Sigma} d| = 1$ and $\mathcal{H}^1(\Sigma \cap \partial C_t)=\mathcal{H}^1(H_C \cap \partial C_t)$, which is only possible when $\Sigma \cap \partial C_t$ are helices, i.e. $\Sigma=H_C$.
\end{proof}

Note that if our competitor surfaces are minimal, transversality condition can be dropped. However, the condition that $\Sigma$ contains the central axis of $C$ cannot be weakened. Without that condition, $\mathcal{H}^1(\Sigma \cap \partial C_t) \geq  \mathcal{H}^1(H_C \cap \partial C_t)$ does not hold for each $t \in(0,r]$ in general; Figure 1(b) can be a counterexample. But if we have a certain condition on the boundary curves of $\Sigma$ instead, we get the same area-minimization result:

\begin{theorem} \label{thm2}
Let $\Sigma \subset C$ be an immersed disk-type surface bounded by a Jordan curve $\Gamma=d_1 \cup d_2 \cup \gamma_1 \cup \gamma_2$ where $\gamma_1$ and $\gamma_2$ are rotationally symmetric $C^2$ curves on S. If the total curvature of $\gamma_1$ is less than $\pi$, then $\mathcal{H}^2 (\Sigma) \geq \mathcal{H}^2 (H_C)$. Equality holds if and only if $\Sigma$ is congruent to $H_C$.
\end{theorem}

\begin{proof}  
Let $\Sigma_0$ be the Douglas-Rad\'o solution for \textcolor{black}{$\Gamma$}. 
Then, $\int_{\Gamma} \kappa  < 4\pi$ implies that $\Sigma_0$ is a unique disk-type minimal surface $\Sigma_0$ by the generalization of Nitsche's uniqueness theorem to curved polygons \cite{Zheng:1990}. Moreover, $\Sigma_0$ is embedded by \cite{Ekholm_White_Wienholtz:2002}.

We claim that $\Sigma_0$ is \textcolor{black}{rotationally symmetric} \textit{i.e.} $\rho_\pi(\Sigma_0)=\Sigma_0$, where $\rho(\cdot)$ is the rotation about the central axis $l$ of $C$ by $\pi$. To show this, suppose $\rho_\pi(\Sigma_0) \neq \Sigma_0$. Then, $\rho_\pi(\Sigma_0)$ and $\Sigma_0$ are two different minimal disks that share the same boundary because $\rho_\pi(\Gamma)=\Gamma$. This contradicts the fact that $\Gamma$ bounds only one minimal disk. Therefore, $\rho_\pi(\Sigma_0)=\Sigma_0$. Moreover, the embeddedness and the \textcolor{black}{rotational} symmetry imply that $\Sigma_0$ should contain the axis of $C$. 

\textcolor{black}{We now show that $\Sigma_0$ is transversal to $\partial C_t=\{(x,y,z) \mid x^2+y^2 \leq t^2\}$ for every $0<t\leq r$. This transversality follows from the fact that a minimal surface cannot touch a cylinder from inside, as follows.  Clearly $\Sigma_0$ is transversal to $\partial C_t$ for sufficiently small $t>0$. Also $\Sigma_0$ is transversal to $\partial C_t$ for $t$ sufficiently close to $r$ because $\Sigma_0\cap \partial C_r$ is a $C^2$ curve. Let $\Sigma_0^1$ be one of the two components of $\Sigma_0\setminus l$. If $\Sigma_0^1\cap \partial C_t$ is a connected curve on $\Sigma_0^1$ for every $0<t\leq r$, then $\Sigma_0^1$ is transversal to $\partial C_t$ for all $0<t\leq r$ and $\Sigma_0^1\cap\partial C_t$ is an analytic curve for every $0<t<r$. This is because if $\Sigma_0^1$ is tangent to $\partial C_t$ at a point $p\in\Sigma_0^1\cap\partial C_t$ then in a neighborhood of $p$, $\Sigma_0^1\cap\partial C_t$ is the union of at least two curves meeting at $p$, and then $\Sigma_0^1\cap\partial C_{t\pm\delta}$ cannot be connected for small $\delta>0$. Suppose $\Sigma_0^1\cap\partial C_t$ is not connected for some $0<t<r$. Then there exists $0<a<r$ such that
	$$a=\inf\{b\,|\,\Sigma_0^1\cap\partial C_t\,\,{\rm is\,\,connected\,\,for\,\,every\,\,} t\in (b,r)\}.$$
	Obviously $\Sigma_0^1\cap\partial C_a$ is not connected. So it consists of at least two components $A$, $B$ such that $E:=\Sigma_0^1\,\cap\, \{(x,y,z)\,|\,a^2< x^2+y^2\leq r^2\}$ connects $A$ and $\partial C_r$ while $B$ is disjoint from $E$. In fact, as $t$ decreases from $r$ to $a$, any point $p\in B$ must be a point of first touching between $\Sigma_0^1\setminus \overline{E}$ and $\partial C_t$. Hence $\Sigma_0^1\setminus \overline{E}$ is tangent to $\partial C_a$ at $p$ from inside, that is, $\Sigma_0^1\setminus\overline{E}\subset C_a$. But this is not possible by the maximum principle. Therefore $\Sigma_0^1\cap\partial C_t$ is a connected and analytic curve for every $0<t<r$ and thus $\Sigma_0$ is transversal to $\partial C_t$ for all $0<t\leq r$.} Therefore, by Theorem \ref{thm1}, we get $\mathcal{H}^2 (\Sigma) \geq \mathcal{H}^2 (\Sigma_0) \geq \mathcal{H}^2 (H_C)$, and it is clear that the equality holds if and only if $\Sigma=H_C$.

\end{proof}

\begin{figure}
\centering
\subfigure[]{\includegraphics[width=0.3\textwidth]{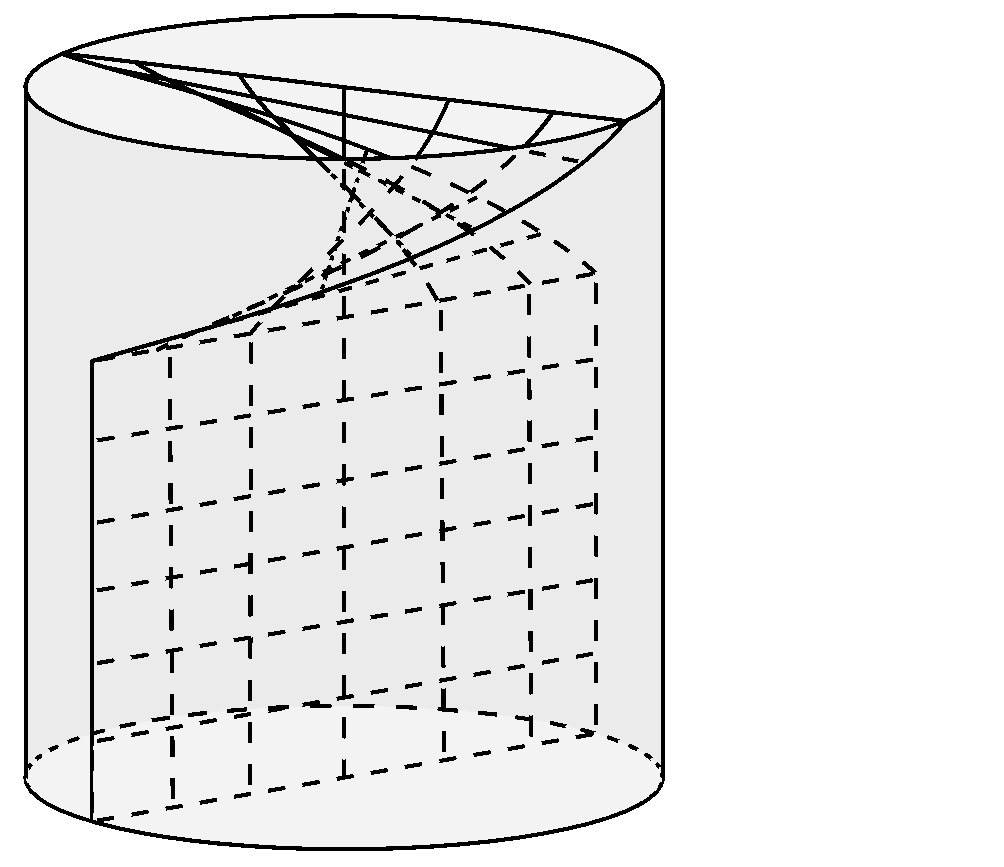}} 
\subfigure[]{\includegraphics[width=0.19\textwidth]{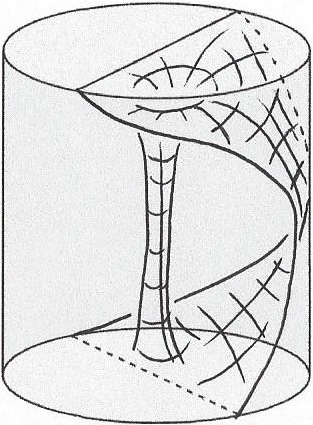}}
\qquad \quad
\subfigure[]{\includegraphics[width=0.21\textwidth]{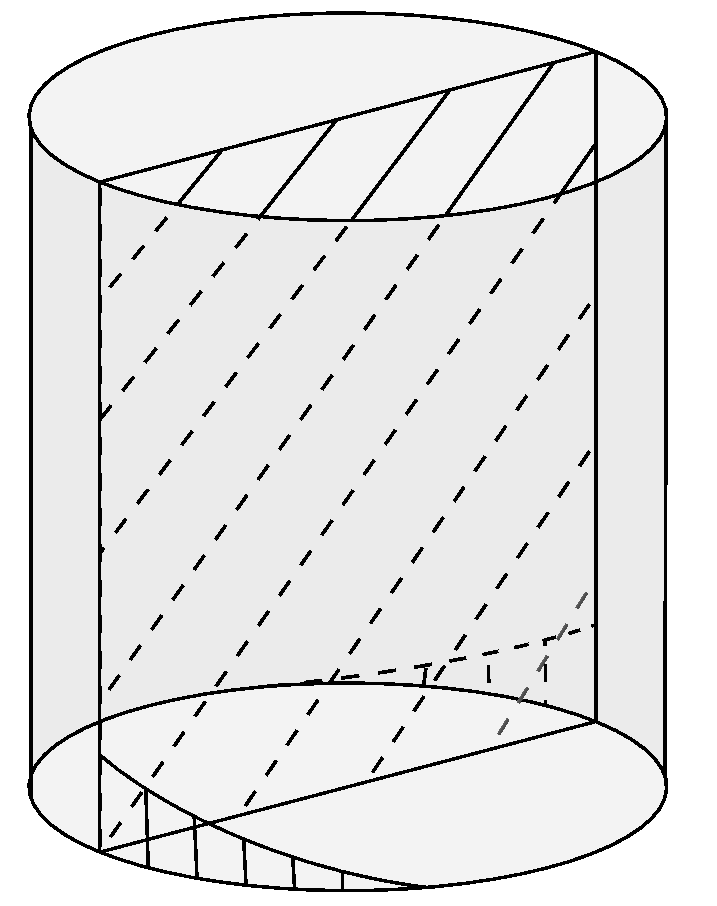}}
\caption{(a) Plane attached to a small helicoidal region; (b) Genus-one surface containing the central axis; (c) Plane with two small triangular regions attached on each side}
\end{figure}

It is noteworthy that, according to Theorem \ref{thm1}, the area of $H_C$ is smaller than the surface composed of a plane rectangle with a small piece of helicoid attached (Figure 2(a)), although most part of it is planar. Theorem \ref{thm1} is also strong in the sense that it holds regardless of the height and radius of the given cylinder $C$. \

Note that the topological condition on $\Sigma$ being a disk-type surface is crucial. In fact, the surface with one handle depicted in Figure 2(b) may have a smaller area than a disk-type surface with the same boundary, even though it contains the central axis of $C$. \

We close this section by mentioning the following two facts. First, in both theorems, the boundary condition that $\Sigma$ contains $d_1$ and $d_2$ as its boundaries cannot be weakened. In fact, if we allow $\Sigma$ to contain only one of the diameters, the surface consisting of a piece of a plane with small triangular regions on $S$ attached on each side (see Figure 2(c)) has a smaller area than $H_C$.

Second, Lawlor \cite{Lawlor:1998} showed a similar area-minimizing property for a compact portion of the half helicoid in a slightly different class of surfaces. Namely, he proved that half of the helicoid centered at the $z$-axis in $\mathbb{R}^3$ has the smallest area among all oriented surfaces \textcolor{black}{that share the same boundary}. From his result, we see that, by the reflection principle of minimal surfaces \cite{Dierkes_Hildebrandt_Sauvigny:2010}, a minimal surface containing the axis of a cylinder and bounded by two pair of helices on the side of cylinder has area greater than \textcolor{black}{or equal to} $H_C$. Theorem \ref{thm1} and Theorem \ref{thm2} not only coincide with this result, but generalize it into a broader class of surfaces: $H_C$ is still the area-minimizer among the surfaces containing the central axis with their boundaries not necessarily helical.

\section{Helicoid as the unique minimal surface}

We now turn our attention to the uniqueness problem for minimal surfaces in $C$ spanning $d_1 \cup d_2 \cup S$. In general, by Douglas' Existence Theorem (Douglas' solution) for the Plateau problem, every Jordan curve in $d_1 \cup d_2 \cup S$ bounds at least one minimal disk. Since there are infinitely many such curves lying on $S$, there exist infinitely many minimal surfaces bounded by $d_1 \cup d_2 \cup S$. But for partially free boundary solutions, meaning that $d_1$ and $d_2$ are the fixed boundaries and the other two boundary curves are free on $S$, $H_C$ might be the unique suface under some conditions. It is clear that $H_C$ satisfies the free boundary condition since it meets $S$ orthogonally along the boundary helices. It is intringuing to investigate when $H_C$ is a unique partially free boundary minimal surface. 

In fact, with no further conditions, there is another minimal surface that meets $S$ orthogonally; there should exist a (locally area-maximizing) minimal surface between $H_C$ and the union $U$ of the top and bottom half disks of $C$ in the one-parameter family of minimal surfaces that are orthogonal to $S$. This is because $H_C$ and $U$ are the relative weak minima in the class of all neighboring surfaces with the same boundary. We will look for some natural boundary conditions under which no partially free boundary solutions of minimal surfaces other than $H_C$ exist. 

Next, we turn out attention to the minimal surfaces in $C$ spanning a pair of helices and will prove that $H_C$ is unique if the angle of the projection of the boundary helices is bounded by the ratio of height and radius of $C$. Another uniqueness theorem follows with no assumption on height and radius of $C$. Note that these two results are not about partially free boundary problem. And, unlike in the previous section, the class of surfaces we consider here are minimal surfaces.

To make our discussion more accurate and general, we need the following definition concerning the curves on the side $S$ of the cylinder $C$. 

\begin{definition} Define the rotation angle $rot(\gamma)$ of a curve $\gamma \subset S$ connecting two points $p, q \in S$ as the oriented angle of the circular arc $proj(\gamma)$ joining $proj(p)$ and $proj(q)$, where $proj(\cdot)$ is the vertical projection of $S$ onto $\partial D_2$.
\end{definition}

\begin{figure}[h]
\vspace{-1cm}
\centering
\includegraphics[height=7.5cm]{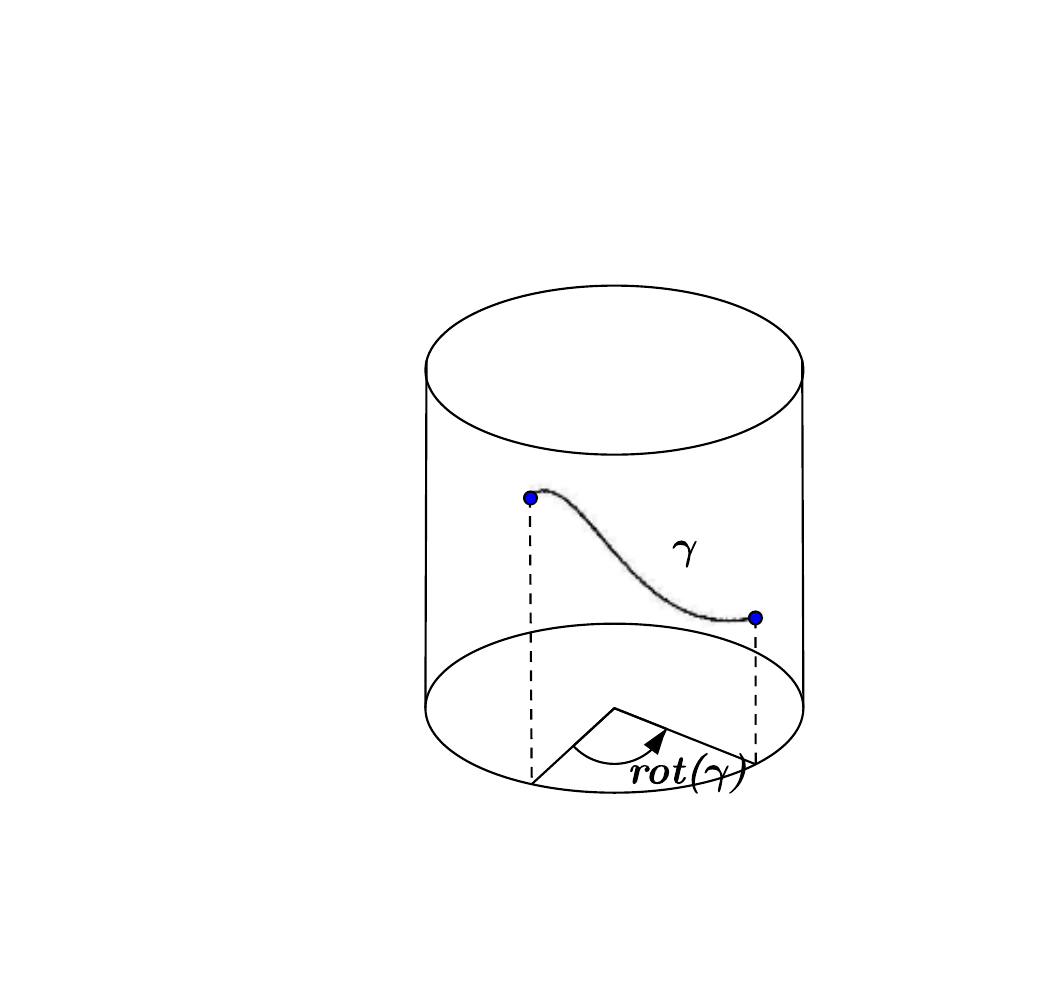}
\vspace{-0.8cm}
\caption {Rotation Angle}
\end{figure}

For example, if a curve $\gamma \subset S$ is part of a helix parameterized by $(r \cos \theta, r \sin \theta, b \, \theta)$ for some constants $r$ and $b$ and $0 \leq \theta \leq \theta_0$, then $rot(\gamma)=\theta_0$. In other words, the rotation angle of a helix measures how much the helix rotates as it sweeps from $d_1$ to $d_2$.

We denote by $H_C^\theta$ the surface $H_C$ with $rot(h_i)=\theta$ where $h_i \in H_C \cap S$ (i=1, 2) are the boundary helices. In fact, $H_C$ in the previous section is actually $H_C^\theta$ with $0 < \theta \leq \frac{\pi}{2}$ or, equivalently, $-\frac{\pi}{2}\leq \theta < 0$. The notion of the rotation angle of the boundary curves was not needed in the previous section because the area of $H_C$ is clearly smaller than that of $H_C^\theta$ with $\theta > \frac{\pi}{2}$ or $\theta < -\frac{\pi}{2}$, so that the area minimizer is $H_C$. 

However, when we investigate the uniqueness of $H_C$, this notation makes it possible to show that not only $H_C$ but the helicoids with boundary helices spiraling up more than $\frac{\pi}{2}$ have a uniqueness property; we shall see that $H_C^\theta$ is the only the minimal surface in $C$ with the rotation angle of the boundary curves on $S$ being $\theta$ under some conditions. Note that thoughout this section we allow $\theta=0$, so that $d_1$ and $d_2$ can be parallel and $H_C^0$ is a plane.

\subsection{$H_C^\theta$ is the Unique Minimal Surface Having Symmetric Free Boundary}

Hereafter, let $\Sigma \subset C$ be a minimal surface spanning $d_1$, $d_2$, and consider two $C^2$ curves $\gamma_1$ and $\gamma_2$ on $S$ such that $\gamma_1 \cap \gamma_2 = \emptyset$ and $\partial (d_1 \cup d_2)=\partial (\gamma_1 \cup \gamma_2)$. With this setting, the question we originally posed was the following: if $\gamma_1, \gamma_2$ are symmetric with respect to the axis $\ell$ of $C$, is $\Sigma$ also symmetric? This remains open, but the next theorem suggests not only an affirmative answer to this question, but leads to a stronger conclusion when $\Sigma$ meets $S$ at a right angle along $\gamma_1$ and $\gamma_2$.

\begin{theorem}\label{thm3}
Let $\Sigma$ be a minimal disk in $C$ with $\partial \Sigma = d_1 \cup d_2 \cup \gamma_1 \cup \gamma_2$, where $\gamma_i \subset S$ are disjoint $C^2$ curves with $rot(\gamma_i)=\theta$ (i=1, 2) and $\Sigma$ is orthogonal to $S$. If $\gamma_1$ and $\gamma_2$ are  \textcolor{black}{rotationally symmetric}, then $\Sigma=H_C^\theta$.
\end{theorem}

The main idea of the proof is based on the maximum principle, which states that two different minimal surfaces cannot touch each other at an isolated interior point. If one minimal surface lies on one side of the other in the neighborhood of that point, the two surfaces must coincide. The same holds for a contact point on the boundary, provided the boundary curve is at least $C^2$. More specifically, the following lemma explains how two different minimal surfaces behave in the neighborhood of their contact points. 

\begin{lemma} \label{lem1}
Let $M$ and $N$ be two minimal surfaces in $\mathbb{R}^3$ that, at a common point $p$, have the same tangent plane $P$. Then, $M$ and $N$ either coincide, or the orthogonal projection of $M \cap N$ on $P$ forms $k$ curved rays emitting from $p$, making the same angles $2\pi/k$ for some even number $k \geq 4$. 
\end{lemma}

\begin{proof} 
See \cite{Fang:1996}.
\end{proof}

\begin{proof}[Proof of Theorem \textcolor{black}{\ref{thm3}}]
Suppose $\Sigma \neq H_C^\theta$. Consider  $\mathcal{F}=\{\rho_\alpha(H_C^\theta):\rho_\alpha$ to be a counter-clockwise rotation function about the axis $\ell$ of $C$ by an angle $\alpha\in[0,\pi)\}$. Then, because $\mathcal{F}$ foliates $C \sim \ell  :=  C \setminus \ell$, $\mathcal{F}$ induces a natural foliation on $\Sigma \sim \ell$ defined by $\mathcal{F} \mid_\Sigma=\{\rho_\alpha(H_C^\theta) \cap \Sigma : \alpha \in [0, \pi)\}$. In other words, $\Sigma$ is filled with the disjoint curves (the leaves of the foliation) $\rho_\alpha(H_C^\theta) \cap \Sigma$. We will verify that these curves and $\Sigma$ itself should satisfy the following four properties:\\

\noindent 
i) $\rho_0(H_C^\theta) \cap (int \Sigma)$ should contain a curve $\gamma_0$ connecting $d_1$ to $d_2$. \\
ii) $\alpha \in (0, \pi)$ $\Rightarrow$ $\rho_\alpha(H_C^\theta) \cap d_i = \{ c_i \}$,  $c_i \phantom{a} \textrm{is the center of} \phantom{a} d_i$. $(i=1,2)$\\
iii) $\alpha_1 \neq \alpha_2$ $\Rightarrow$ $\rho_{\alpha_1} (H_C^\theta) \cap \Sigma$ and  $\rho_{\alpha_2}(H_C^\theta) \cap \Sigma$ are disjoint except on $\ell$. \\
iv) $\exists$ $\alpha_0 \in (0,\pi)$ such that
\begin{displaymath}
\left \{ \begin{array}{ll}
0 \leq \alpha \leq \alpha_0 \Rightarrow & \rho_\alpha(H_C^\theta) \cap \gamma_i \neq \emptyset \ \textrm{and $\rho_{\alpha_0} (H_C^\theta)$ is tangent to $\Sigma$} \\
& \textrm{at some point $p_i \in \gamma_i$.} \\
\alpha > \alpha_0 \phantom{>>}\Rightarrow & \rho_\alpha(H_C^\theta) \cap \gamma_i = \emptyset.
\end{array} \right.
\end{displaymath}

Property i) follows from the boundary conditions that $d_1 \cup d_2 \subset \Sigma \cap H_C^\theta$ and the $\gamma_i$ are  \textcolor{black}{rotationally symmetric} with each other. Properties ii) and iii) are straightforward because $\Sigma \sim \ell$ is foliated by $\{\rho_\alpha(H_C^\theta) \cap \Sigma : \alpha \in [0, \pi)\}$. To see iv), we proceed as follows: \\
Let $\partial d_i=\{p^i, q^i\}$ $(i=1,2)$. Develop $S$ onto $\mathbb{R}^2$ after cutting it at the vertical plane passing through $d_1$. Then, we have a rectangular region with points $p^i$ and $q^i$ on the top and bottom edges. Since $rot(\gamma_i)=rot(h_i)=\theta$, the developed curves of the $\gamma_i$ and $h_i$ have the same end points, connecting the top edge to the bottom. Note that the developed curves of $h_i$ are actually line segments on the rectangular region. Hence, by the Mean Value Theorem, there exists at least one point, say $p$, on $\gamma_1$ and another point $p'$ on $\gamma_2$ at which the tangent lines have the same slope as the $h_i$. Since $p$ and $p'$ are the common points of $\Sigma$ and $\rho_\alpha(H_C^\theta)$, and since $\Sigma$ meets  \textcolor{black}{$S$} orthogonally, $p$ and $p'$ become contact points of $\Sigma$ and $\rho_\alpha(H_C^\theta)$ at which they have the same tangent plane. Moreover, there will be no more intersection points on $\gamma$ if we choose $\alpha_0$ to be the maximum angle corresponding to the last contact point. This proves iv). 

We will show that property iv) and the maximum principle for minimal surfaces imply the leaves $\rho_\alpha(H^{\theta}_C) \cap \Sigma$ of the foliation of $\Sigma$ should behave in a special way, leading to a contradiction. Let us first look at the case of only one pair of boundary contact points $p_1 \in \gamma_1$, $p_2 \in \gamma_2$ of $\rho_{\alpha} (H_C^\theta)$ and $\Sigma$ at $\alpha=\alpha_0$. 
By ii) (that is, $\rho_\alpha(H_C^\theta) \cap d_i$ is a single point), there should be exactly one curve emitting from $d_i$ for each $i=1, 2$. Because of the axial symmetry of the $\gamma_i$, $\rho_\alpha(H_C^\theta) \cap \partial \Sigma$ are pairwise symmetric points on $\gamma_1$ and $\gamma_2$ with respect to the same $\alpha$, converging to $p_1$ and $p_2$ as $\alpha$ goes to $\alpha_0$.

By the assumption that $\Sigma$ is not congruent to $H_C^\theta$ and property iv), Lemma \ref{lem1} implies that near $p_1$, $\rho_{\alpha_0} (H^\theta_C) \cap \Sigma$ consists of $k$ curves emitting from $p_1$, making an angle of $\pi/k$ ($k \geq 2$). We assume, for now, that $k=2$.  


Note that the leaves $\rho_\alpha(H_C^\theta) \cap \Sigma$ should not be a closed curve in $\Sigma \sim \partial \Sigma$ for each $\alpha$. Instead, every leaf should be connected to the boundary curves $\gamma_i$. This is because, if there is a region enclosed by the intersection curve, then those intersection curves should converge to an interior point as $\alpha$ increases. On that point, we can apply the maximum principle to get a contradiction.

It is now obvious that, under i), iii), and the previous observations, there exists $\alpha_C \in (\alpha_0, \pi)$ such that $\rho_{\alpha_C}(H^\theta_C) \cap int(\Sigma)$ is a single curve, say $\gamma_C$, along which $\rho_{\alpha_C}(H^\theta_C)$ and $\Sigma$ are tangent to each other. In other words, $\rho_{\alpha}(H^\theta_C) \cap int(\Sigma)$ converge to $\gamma_C$, which is one leaf of the foliation, such that $\rho_{\alpha_C}(H^\theta_C)$ is located on one side of $\Sigma$ along $\gamma_C$. However, this contradicts the boundary maximum principle. Therefore, we conclude that $\Sigma$ must be equal to $H^\theta_C$ when $k=2$. 

If $k \geq 3$, we have the same conclusion, as the same argument can be applied to the component that contains the leaves $\rho_{\alpha} (H^\theta_C) \cap \Sigma$ of $\alpha > \alpha_0$. Specifically, near $p_1$, $\Sigma$ is divided by $k+1$ regions such that, if one component is comprised of $\rho_{\alpha} (H^\theta_C) \cap \Sigma$ of $\alpha < \alpha_0$, then the adjacent component is of $\rho_{\alpha} (H^\theta_C) \cap \Sigma$ of $\alpha > \alpha_0$. The behavior of the leaves of $\rho_{\alpha} (H^\theta_C) \cap \Sigma$ is the same as for $k=2$, thus eliciting the same contradiction. In other words, we have $\Sigma=H^{\theta}_C$ when there are two boundary contact points of $\rho_{\alpha_0} (H^\theta_C)$ and $\Sigma$ (one at each $\gamma_i$).

In general cases, that is, when the number of the boundary contact points is greater than two, we get the same contradiction in the region containing $p_i$, the last boundary contact points of $\rho_{\alpha} (H^\theta_C)$ and $\Sigma$. 
\end{proof}

$Remark.$ From the proof, it is obvious that when $d_1$ and $d_2$ are parallel with $rot(\gamma_i)=0$, $\Sigma$ should be a piece of the plane under the same hypothesis. 
Also, note that the minimality of $\Sigma$ is crucial in the proof: if $\Sigma$ is not minimal, the intersection curves can be arbitrary, and one cannot assert the existence of the converging curve along which $\rho_{\alpha_C}(H^\theta_C)$ is located on one side of $\Sigma$.

\subsection{$H_C^\theta$ is the Unique Surface Spanning Helices}

We now turn to the case in which $\partial \Sigma \cap S$ is a double helix and $\Sigma$ does not necessarily meet $S$ orthogonally. A double helix is a pair of helices that are symmetric to each other with respect to the axis of $C$. For simplicity, but without losing the generality, let us assume that the rotation angle $\theta$ is non-negative in this section. Then the next theorem asserts that such $\Sigma$ should be equal to part of the helicoid if the height $h$ is bigger than the product of the radius $r$ and the rotation angle $\theta$.

\begin{theorem}\label{thm4}
Suppose $\Sigma \subset C$ is a minimal surface spanning $d_1, d_2$ and a double helix $h_1, h_2$ with $rot(h_i)=\theta$. If $0 \leq \theta \leq h/r$, then $\Sigma = H_C^\theta$. 
\end{theorem}

\begin{proof} Assume that $\Sigma \neq H_C^\theta$. Let $\Omega \in C$ be the surface obtained by the screw motion (along the $z$-axis) of  \textcolor{black}{a circular} arc (in the $xy$-plane) with radius $R>0$ and central angle $2\eta \in (0, \pi]$. Then, for $a:=h/\theta$, $\Omega$ can be parametrized by 
\begin{align*}
\Omega (\zeta, \phi)&=(R \cos \eta \sin \zeta-R \cos(\pi/2-\zeta+\eta-\phi), -R \cos \eta \cos \zeta+R \sin(\pi/2-\zeta+\eta-\phi), a\zeta)\\
&=(-R \sin (\zeta-\eta+\phi)+R \cos \eta \sin \zeta, R \cos (\zeta-\eta+\phi)-R \cos \eta \cos\zeta, a\zeta),
\end{align*}
where $\sin \eta = \frac{r}{R}$. See Figure \ref{fig:geo}. \\ 

\noindent A computation yields
\begin{align*}
&\Omega_\zeta=(-R\cos(\zeta-\eta+\phi)+R\cos \eta \cos\zeta, -R\sin(\zeta-\eta+\phi)+R\cos\eta\sin\eta, a) \\
&\Omega_\phi=(-R\cos(\zeta-\eta+\phi), -R\sin(\zeta-\eta+\phi), 0) \\
&\Omega_\zeta \times \Omega_\phi=(aR\sin(\zeta-\eta+\phi),-aR\cos(\zeta-\eta+\phi),R^2\cos\eta\sin(\eta-\phi)) \\
&\vec{N}=\frac{(a \sin (\zeta-\eta+\phi), -a \cos(\zeta-\eta+\phi), R \cos \eta \sin (\eta-\phi))}{\sqrt{a^2+\rho^2 \cos^2 \eta \sin^2 (\eta-\phi)}} \\
&\Omega_{\zeta \zeta}=(R\sin(\zeta-\eta+\phi)-R\cos\eta\sin\zeta, -R\cos(\zeta-\eta+\phi)+R\cos\eta\cos\zeta, 0)\\
&\Omega_{\zeta \phi}=(R\sin(\zeta-\eta+\phi), -R\cos(\zeta-\eta+\phi),0)=\Omega_{\phi \phi}.
\end{align*} 
Therefore, the coefficients of the first and second fundamental forms of $\Omega$ are: 
\begin{align*}
g_{11}=<\Omega_\zeta, \Omega_\zeta>&=R^2(1+\cos^2\eta-2\cos\eta\cos(\eta-\phi))+a^2 \\
g_{12}=<\Omega_\zeta, \Omega_\phi>&=R^2-R^2\cos\eta\cos(\eta-\phi) \\
g_{22}=<\Omega_\phi, \Omega_\phi>&=R^2\\
b_{11}=<\Omega_{\zeta\zeta}, \vec{N}>&=\frac{aR(1-\cos\eta\cos(\eta-\phi))}{\sqrt{a^2+R^2 \cos^2 \eta \sin^2 (\eta-\phi)}} \\
b_{12}=<\Omega_{\zeta\phi}, \vec{N}>&=\frac{aR}{\sqrt{a^2+R^2 \cos^2 \eta \sin^2 (\eta-\phi)}} \\
b_{22}=<\Omega_{\phi\phi}, \vec{N}>&=\frac{aR}{\sqrt{a^2+R^2 \cos^2 \eta \sin^2 (\eta-\phi)}},
\end{align*}
\begin{align*}
\therefore \textrm{H}(\Omega)&=\frac{1}{2}\frac{g_{22}b_{11}-2g_{12}b_{12}+g_{11}b_{22}}{g_{11}g_{22}-g_{12}^2}\\
&=\frac{a(a^2+R^2\cos^2\eta-R^2\cos \eta \cos (\eta-\phi))}{{2R(a^2+R^2 \cos^2 \eta \sin^2 (\eta-\phi))}^{3/2}}\\
&=\frac{a(a^2-r^2+R^2(1-\cos \eta \cos (\eta-\phi)))}{{2R(a^2+R^2 \cos^2 \eta \sin^2 (\eta-\phi))}^{3/2}}.
\end{align*}
Substituting $\sin \eta = \frac{r}{R}$ induces the last equality. 

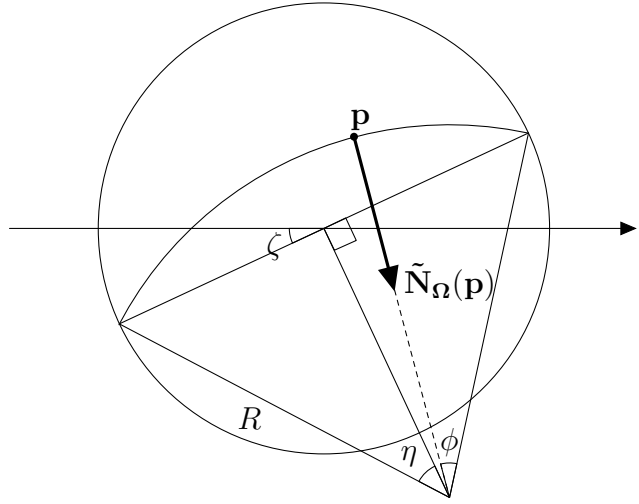
\begin{figure}
\centering
\begin{tikzpicture}[line cap=round,line join=round,>=triangle 45,x=1.0cm,y=1.0cm]
\clip(-4.3,-4.24) rectangle (4.51,4.05);
\draw [shift={(0,0)}] (0,0) -- (180:0.46) arc (180:205.1:0.46) -- cycle;
\draw(0.14,-0.29) -- (0.43,-0.16) -- (0.29,0.14) -- (0,0) -- cycle; 
\draw [shift={(1.67,-3.58)}] (0,0) -- (115.03:0.46) arc (115.03:152.29:0.46) -- cycle;
\draw [shift={(1.67,-3.58)}] (0,0) -- (77.77:0.46) arc (77.77:104.87:0.46) -- cycle;
\draw(0,0) circle (3cm);
\draw (-2.72,-1.27)-- (2.72,1.27);
\draw (0,0)-- (1.67,-3.58);
\draw (-2.72,-1.27)-- (1.67,-3.58);
\draw (2.72,1.27)-- (1.67,-3.58);
\draw [shift={(1.67,-3.58)}] plot[domain=1.36:2.66,variable=\t]({1*4.96*cos(\t r)+0*4.96*sin(\t r)},{0*4.96*cos(\t r)+1*4.96*sin(\t r)});
\draw [dash pattern=on 2pt off 2pt] (0.4,1.22)-- (1.67,-3.58);
\draw [->,line width=1.2pt] (0.4,1.22) -- (0.94,-0.84);
\draw (-0.99, -2.53) node {$R$};
\draw (-0.91,0.1) node[anchor=north west] {$\mathbf{\zeta}$};
\draw (0.88,-2.76) node[anchor=north west] {$\mathbf{\eta}$};
\draw (1.40,-2.51) node[anchor=north west] {$\mathbf{\phi}$};
\draw (0.92,-0.32) node[anchor=north west]{$\mathbf{\vec{N}_\Omega(p)}$};
\draw (0.21,1.72) node[anchor=north west] {$\mathbf{p}$};
\draw [->] (-4.18,0) -- (4.16,0);
\begin{scriptsize}
\fill [color=black] (0.4,1.22) circle (1.5pt);
\end{scriptsize}
\end{tikzpicture}
\caption[Parameterization of a Barrier Surface]{Parameterization of $\Omega$ at height $a\zeta$} \label{fig:geo}
\end{figure}

Observe that the mean curvature $\textrm{H}(\Omega)$ of $\Omega$ is always greater than zero when $h\geq r\theta$, that is, $a=\frac{h}{\theta} \geq r$. Since the normal vector $\vec{N_\Omega}$ of $\Omega$ is always pointing toward $\Sigma$ (see Figure \ref{fig:geo}), so is the mean curvature vector $\vec{\textrm{H}}(\Omega)=\textrm{H}(\Omega)\vec{N_\Omega}$ at any point of $\Omega$. Therefore, if there is a contact point of $\Omega$ and $\Sigma$ inside $C$, we can apply the interior maximum principle there to get a contradiction. Fortunately, this is possible by choosing $\eta$ appropriately and, if necessary, taking the surface $\Omega '$ obtained by reflecting $\Omega$ with respect to $H_C^\theta$ instead of $\Omega$; as $\eta$ varies, $\Sigma$ must have a contact point with either $\Omega$ or $\Omega '$ for some $\eta \in (0, \pi/2)$. This contradicts the maximum principle. In other words, $\Omega = H_C^\theta$ if $h \geq r\theta$. 
\end{proof}

 \textcolor{black}{One important remark is that our proof includes wider classes of surfaces. They do not have to be topological disks. Also, their boundary helices can have the total curvature bigger than $2\pi$ as long as $\theta \leq h/r$. In fact, if the total curvature of the boundary helices $h_1, h_2$ \textcolor{black}{is less than $2\pi$}, the total curvature of $d_1 \cup d_2 \cup h_1 \cup h_2$ \textcolor{black}{is less than} $4\pi$, which implies that the disk-type minimal surface bounded by $d_1 \cup d_2 \cup h_1 \cup h_2$ is unique by \textcolor{black}{a generalization of Nitsche's uniqueness theorem \cite{Zheng:1990}} and therefore should be part of the helicoid.}  

If there is no restriction on the ratio of the height and radius of $C$, the method used in the proof of Theorem \ref{thm4} can no longer be applied, because the mean curvature of $\Omega$ may be negative. However, we obtain the same uniqueness result with no conditions on the height and radius of $C$. Namely, if $0 \leq rot(h_i) \leq \pi$, $\Sigma$ should be congruent to $H_C^\theta$.

\begin{theorem}\label{thm5}
If $\Sigma \subset C$ is a disk-type minimal surface spanning $d_1$, $d_2$, $h_1$, and $h_2$ with $0 \leq rot(h_i) \leq \pi$, then $\Sigma = H_C^\theta$ with $0 \leq \theta \leq \pi$.
\end{theorem}

\begin{proof}  \textcolor{black}{Let $rot(h_i)=\phi_0$. Without loss of generality, suppose $\textcolor{black}{d_2} (t)=(rt, 0, 0)$, $\textcolor{black}{d_1} (t)=(rt\cos \phi_0, rt\sin \phi_0, b\phi_0)$, $h_1(u)=(r\cos u, r\sin u, bu)$, and $h_2(u)=(-r\cos u, -r\sin u, bu)$  for $-1 \leq t \leq 1$, $0 \leq u \leq \phi_0$ and for some positive constant $b$. Let $\pi_{yz}$ be the projection map onto the $yz$-plane. Then, $\pi_{yz} \circ h_1 (u)=(0, r\sin u, bu)$, $\pi_{yz} \circ h_2 (u)=(0, -r\sin u, bu)$, $\pi_{yz} \circ d_1 (t)=(0, 0)$, and $\pi_{yz} \circ d_2 (t)=( rt\sin \phi_0, b\phi_0)$. Hence, $rot(h_i)=\phi_0 \leq \pi$ implies that $\partial \Sigma=d_1 \cup d_2 \cup h_1 \cup h_2$ is a Jordan curve which has a monotonic orthogonal projection onto a convex plane Jordan curve. Then, by the Remark on p299 of \cite{Dierkes_Hildebrandt_Sauvigny:2010} or by Theorem 2 of \cite{Meeks:1981} (a generalized version of Rado's Theorem), $\partial \Sigma$ bounds a unique minimal disk. To be more specific, because Rado's Theorem remains true when vertical segments of a given Jordan curve are mapped onto single points of a plane convex curve by an orthogonal projection, we conclude that $H_C^\theta$ is the unique minimal disk spanning $d_1$, $d_2$, $h_1$, and $h_2$ with $0 \leq rot(h_i) \leq \pi$.}    
\end{proof}

Note that the method used in the proof of Theorem \ref{thm5} is not valid when $rot(h_i)>\pi$, because the region bounded by the projection map is no longer convex. Theorem \ref{thm5} is related to Theorem \ref{thm4} in the sense that both provide an upper bound of $\theta$ to assert the uniqueness of $H_C^\theta$. The differences are that Theorem \ref{thm4} uses the ratio of $h$ over $r$ to bound $\theta$ with no topological restrictions, whereas Theorem \ref{thm5} reaches the same conclusion regardless of $h$, $r$, or $h/r$ with a topological restriction. \\

\subsection{Remarks}

In \cite{Nitsche:1989}(\S111), Schwarz's result on the stable part of the helicoid is introduced. Namely, the area of the helicoid with pitch $2\pi$: $\{(u\cos v, u\sin v, v): -r \leq u \leq r, -\frac{\theta}{2} \leq v \leq \frac{\theta}{2} \}$ is a relatively weak minimum in the class of all neighboring surfaces with the same boundary if $\theta \leq \pi$ or $\theta > \pi$ and $r \leq r(\theta)$. The upper bound of $r(\theta)$ is given by $\bar{r} \approx 1.5088 \cdots$. Beware that the range of $\theta$ there is not exactly the same as ours. \

Theorem \ref{thm4} differs from this result in the sense that Schwarz's result was about the stability, whereas our focus is on uniqueness. To be more precise, a certain piece of the helicoid has, according to Schwarz's result, a smaller area than the {\it neighboring} surfaces that have helices as their boundary, whereas Theorem \ref{thm4} guarantees that this piece of the helicoid is, in fact, unique among {\it all} minimal surfaces with the same boundary. In other words, Theorem \ref{thm4} extends Schwarz’s result by asserting that, when $r \leq 1$, the helicoid with pitch $2\pi$ is not only stable but also unique. See the shaded region in Figure 5(a). Although $r$ is not covered up to $\bar{r} \approx 1.5088 \cdots$, our proof is more geometrical and avoids long calculations. Also, our result is independent of the pitch of th helicoid. We possibly reach $\bar{r}$ if we can obtain the lower bound of the mean curvature $\textrm{H}(\Omega)$ in the proof of Theorem \ref{thm4}. \\

\begin{figure}\label{graph}
\centering
\subfigure[]{\includegraphics[width=0.5\textwidth]{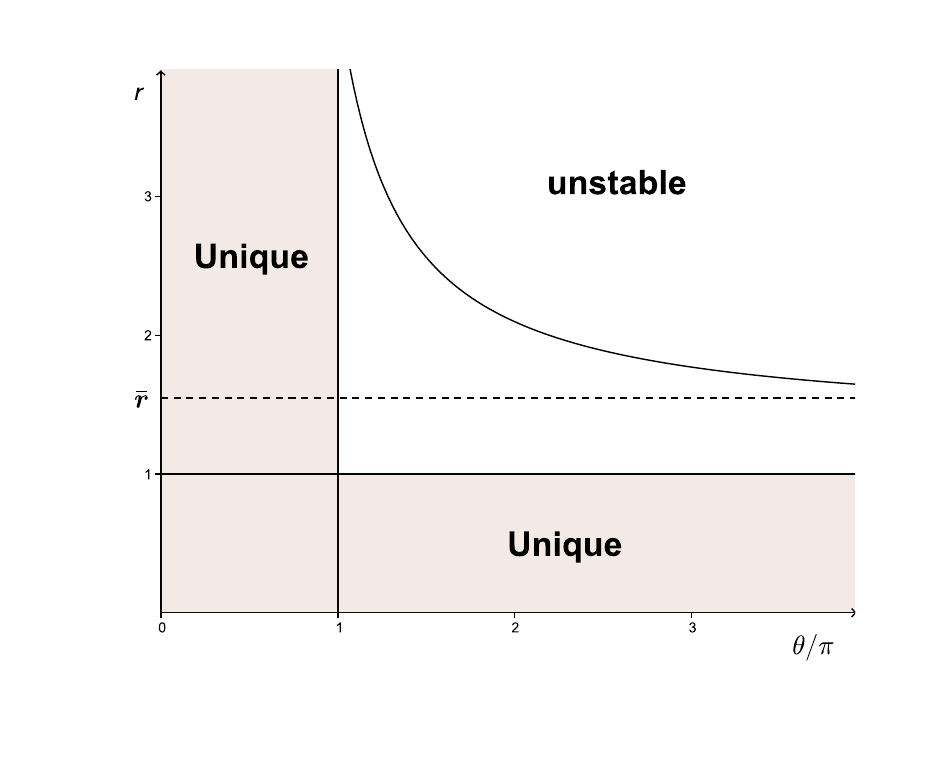}}
\quad
\subfigure[]{\includegraphics[width=0.4\textwidth]{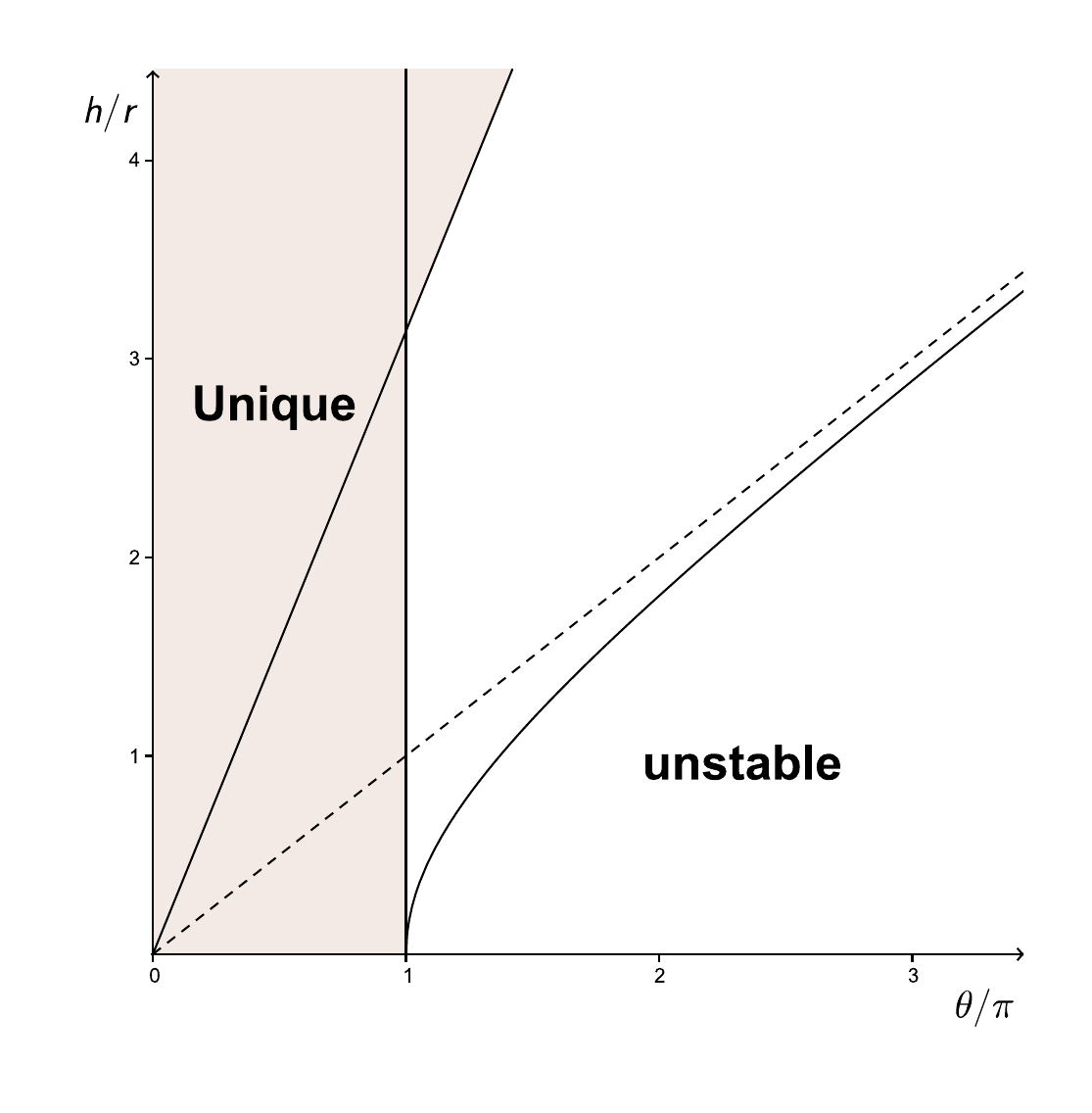}}
\qquad
\caption{(a) The phase diagram of the helicoid of pitch $2\pi$; (b) The phase diagram of the helicoid as the pitch changes}
\end{figure}

In \cite{Boudaoud_Patricio_Amar:1999}, a similar result concerning the stable part of complete helicoid was studied in a slightly different setting: instead of fixing the pitch of the helicoid, the authors obtained the phase diagram of $h/r$ with respect to $\theta$ for which the helicoid was stable, where $h$ and $r$ are the height and radius of the helicoid, respectively. Therefore, Theorem \ref{thm4} is stronger than their result, because it asserts that part of the helicoid is not only stable, but is actually unique when $h/r \geq \theta$ (see Figure 5(b)). \\

After completing the work reported in this paper, the author has learned that \cite{Ripoll_Tomi:1995} contains a similar uniqueness result for the stable part of the helicoid. Namely, it was proved that the stable part of the helicoid is unique among the properly immersed minimal surfaces spanning a double helix, having bounded curvature and satisfying the following asymptotic condition
\[
\lim_{(x,y,z) \to \infty} \frac{x^2+y^2}{z^2}=0.
\]

The main differences are that again we do not fix the pitch of the boundary helices to $2\pi$, nor do we require our surfaces to have bounded curvature or satisfy the asymptotic condition. Our method is also simpler in the sense that we use the surface $\Omega$ obtained by the screw motion of an arc in the plane as a barrier for applying the maximum principle, whereas they used the one-parameter family of associate minimal surfaces of the helicoid. Unlike \cite{Ripoll_Tomi:1995}, $\Sigma$ shares the same boundary as $\Omega$, and it is therefore geometrically clear that an interior touching point exists. In \cite{Ripoll_Tomi:1995}, it was necessary to prove the existence of the interior touching points, as the boundary of their barrier surfaces are not always on one side of $C$. \\

\subsection{Acknowlegement}
The author would like to thank the reviewer for constructive comments and suggestions that improved the manuscript a lot.


\begin{thebibliography}{99}

\bibitem{Bernstein_Breiner:2014} J. Bernstein and C. Breiner, \textit{A variational characterization of the catenoid}, Calc. Var. Partial Differential Equations 49 (2014), 215--232. 

\bibitem{Boudaoud_Patricio_Amar:1999} A. Boudaoud, P. Patricio and M. B. Amar, \textit{The helicoid versus the catenoid: Geometrically Induced Bifurcations}, Phys. Rev. Lett. 83 (1999), 3836--3839. 

\bibitem{Catalan:1842} E. Catalan, \textit{Sur les surfaces r\'egl\'ees dont l'aire est un minimum}, J. Math. Pure Appl. 7 (1842), 203--211.

\bibitem{Chavel:2001} I. Chavel, \textit{Isoperimetric inequalities. Differential geometric and analytic perspectives}, Cambridge Tracts in Mathematics 145, Cambridge University Press, Cambridge, 2001.

\bibitem{Choe_Hoppe:2013} J. Choe and J. Hoppe, \textit{Higher dimensional minimal submanifolds generalizing the catenoid and helicoid}, Tohoku Math. J. (2) 65 (2013), no. 1, 43--55.

\bibitem{Colding_Minicozzi:2004} T.H. Colding and W. P. Minicozzi II, \textit{The space of embedded minimal surfaces of fixed genus in a 3-manifold I, Estimates off the axis for disks}, Ann. of Math. (2) 160 (2004), no. 1, 27--68. 

\bibitem{Collin_Krust:1991} P. Collin and R. Krust, \textit{Le problème de Dirichlet pour l'équation des surfaces minimales sur des domaines non bornés}, Bull. Soc. Math. France 119 (1991), 443--462.


\bibitem{Ekholm_White_Wienholtz:2002} T. Ekholm, B. White and D. Wienholtz, \textit{Embeddedness of minimal surfaces with total boundary curvature at most $4\pi$}, Ann. of math. (2) 155 (2002), no. 1, 209--234.  

\bibitem{Fang:1996} Y. Fang, \textit{Lectures on minimal surfaces in $\mathbb{R}^3$}, Proceedings of the Centre for Mathematics and its Applications, Australian National University 35, Australian National University, Centre for Mathematics and its Applications, Canberra, 1996.

\bibitem{Giusti:1984} E. Giusti, \textit{Minimal surfaces and functions of bounded variation}, Monographs in Mathematics 80, Birkhäuser Verlag, Basel, 1984.

\bibitem{Dierkes_Hildebrandt_Sauvigny:2010} R. Jakob, U. Dierkes, A. K{\"u}ster, S. Hildebrandt and F. Sauvigny, \textit{Minimal Surfaces}, Grundlehren der mathematischen Wissenschaften 295, Springer Berlin Heidelberg, 2010. 

\bibitem{Lawlor:1998} G. Lawlor, \textit{Proving area minimization by directed slicing}, Indiana Univ. Math. J. 47 (1998), no. 4, 1547--1592.

\bibitem{Lee_Lee:2017} E. Lee and H. Lee,  \textit{Generalizations of the Choe–Hoppe helicoid and Clifford cones in Euclidean space}, J. Geom. Anal. 27 (2017), Issue 1, 817--841

\bibitem{Meeks:1981} W. H. Meeks III, \textit{Uniqueness theorems for minimal surfaces}, Illinois J. Math. 25 (1981), no. 2, 318--336. 

\bibitem{Meeks_Rosenberg:2005} W. H. Meeks III and H. Rosenberg, \textit{The uniqueness of the helicoid}, Ann. of Math. (2) 161 (2005), no.2, 727--758.

\bibitem{Morgan:2000} F. Morgan, {\it Geometric Measure Theory. A beginner's guide}, 3rd edition. Academic Press, Inc., San Diego, CA, 2000. 



\bibitem{Nitsche:1989} J.C.C. Nitsche, {\it Lectures on Minimal Surfaces}, Cambridge University Press, 1989.

\bibitem{Rado:1930} T. Rado, {\it Some remarks on the problem of Plateau}, Proc. Natl. Acad, Sci. USA 16 (1930), 242--248

\bibitem{Ripoll_Tomi:1995} J. Ripoll and F. Tomi, \textit{Maximum principles for minimal surfaces in $\mathbb{R}^ 3$ having noncompact boundary and a uniqueness theorem for the helicoid}, Manuscripta Math. 87 (1995), no. 4, 417--434.

\bibitem{Zheng:1990} S.Y. Zheng, \textit{$4\pi$ Uniquess theorem for curved polygons}, Analysis 10 (1990), 247--264.


\end{thebibliography}
\end{document}